\documentclass[11pt]{article}
\usepackage{enumerate}
\usepackage{amssymb,a4wide,latexsym,makeidx,epsfig,fleqn}
\usepackage{amsthm}
\usepackage{amsmath}
\usepackage{enumerate}
\usepackage{graphicx}
\usepackage{float}
\usepackage{indentfirst}
\usepackage[numbers,sort&compress]{natbib}
\allowdisplaybreaks[4]
\newtheorem{theorem}{Theorem}[section]

\newtheorem{lemma}[theorem]{Lemma}

\newtheorem{problem}[theorem]{Problem}

\begin{document}
\textwidth 150mm \textheight 225mm
\title{On the $\alpha$-index of minimally 2-connected graphs with given order or size \thanks{Supported by the National Natural Science Foundation of China (No. 12271439).}}
\author{{Jiayu Lou$^{a,b}$, Ligong Wang$^{a,b,}$\thanks{Corresponding author.}, Ming Yuan$^{a}$}\\
{\small $^a$School of Mathematics and Statistics, Northwestern
Polytechnical University,}\\ {\small  Xi'an, Shaanxi 710129,
P.R. China.}\\
{\small $^b$ Xi'an-Budapest Joint Research Center for Combinatorics, Northwestern
Polytechnical University,}\\
{\small Xi'an, Shaanxi 710129,
P.R. China. }\\
{\small E-mail: jyloumath@163.com, lgwangmath@163.com, ym19980508@mail.nwpu.edu.cn} }
\date{}
\maketitle
\begin{center}
\begin{minipage}{120mm}
\vskip 0.3cm
\begin{center}
{\small {\bf Abstract}}
\end{center}
{\small For any real $\alpha \in [0,1]$, Nikiforov defined the $A_\alpha$-matrix of a graph $G$ as $A_\alpha(G)=\alpha D(G)+(1-\alpha)A(G)$, where $A(G)$ and $D(G)$ are the adjacency matrix and the diagonal matrix of vertex degrees of $G$, respectively. The largest eigenvalue of $A_\alpha(G)$ is called the $\alpha$-index or the $A_\alpha$-spectral radius of $G$. A graph is minimally $k$-connected if it is $k$-connected and deleting any arbitrary chosen edge always leaves a graph which is not $k$-connected. In this paper, we characterize the extremal graphs with the maximum $\alpha$-index for $\alpha \in [\frac{1}{2},1)$ among all minimally 2-connected graphs with given order or size, respectively.

 }
\end{minipage}
\end{center}

\section{Introduction }
\label{sec:ch6-introduction}
\noindent Let $G=(V(G),E(G))$ be a simple undirected graph with vertex set $V(G)$ and edge set $E(G)$. Let $n=|V(G)|$ and $m=|E(G)|$ denote the order and size of $G$, respectively. For a vertex $v \in V(G)$, its neighbor set is denoted by $N_G(v)$ (or, $N(v)$ for short) and its closed neighbor set is defined as $N_G[v]=N_G(v)\cup \{v\}$ (or, $N[v]$ for short). The degree of vertex $v$ is denoted by $d_G(v)=|N_G(v)|$ (or, $d(v)$ for short). Let $\Delta(G)$ and $\delta(G)$ be the maximum degree and minimum degree of $G$, respectively. For a vertex set $S \subseteq V(G)$, let $G[S]$ be the subgraph of $G$ induced by $S$. For $A, B \subset V(G)$, we denote by $e(A)$ the number of edges in $G[A]$ and by $e(A, B)$ the number of edges with one endpoint in $A$ and one endpoint in $B$. Let $G-v$ denote the graph obtained from $G$ by deleting the vertex $v$ together with all the edges incident with $v$. Similarly, Let $G-uv$ (resp. $G+uv$) denote the graph obtained from $G$ by deleting (resp. adding) the edge $uv\in E(G)$ (resp. $uv \notin E(G)$). Let $K_{s,t}$ be a complete bipartite graph with bipartition $(X,Y)$, where $|X|=s$ and $|Y|=t$. For an odd number $m$, $SK_{2,\frac{m-1}{2}}$ (see Fig.\;1) denotes the graph obtained from the complete bipartite graph $K_{2,\frac{m-1}{2}}$ by subdividing one edge. A cycle $C$ of $G$ is said to have a chord if there is an edge of $G$ that joins a pair of non-adjacent vertices from $C$.\indent\setlength{\parindent}{1em}

The adjacency matrix $A(G)$ of $G$ is an $n\times n$ matrix $(a_{ij})_{n\times n}$, where $a_{ij}=1$ if $v_iv_j\in E(G)$ and 0 otherwise. Let $D(G)$ be the diagonal matrix of vertex degrees of $G$. The signless Laplacian matrix of $G$ is defined as $Q(G) =D(G)+A(G)$. The largest eigenvalue of $A(G)$ is called the index or the spectral radius of $G$, and the largest eigenvalue of $Q(G)$ is called the $Q$-index or the signless Laplacian spectral radius of $G$. For any real $\alpha \in [0,1]$, Nikiforov \cite{N2017} proposed to study the convex linear combinations $A_\alpha(G)$ of $A(G)$ and $D(G)$ defined by
$$A_\alpha(G)=\alpha D(G)+(1-\alpha)A(G).$$
It is easy to see that $A(G)=A_0(G)$, $D(G)=A_1(G)$ and $Q(G)=2A_{\frac{1}{2}}(G).$ The largest eigenvalue of $A_\alpha(G)$, denoted by $\rho_\alpha(G)$, is called the $\alpha$-index or the $A_\alpha$-spectral radius of $G$. For a connected graph $G$, $A_\alpha(G)$ is irreducible. By the Perron-Frobenius Theorem, $\rho_\alpha(G)$ is positive, and there exists a unique positive unit eigenvector corresponding to $\rho_\alpha(G)$, which is called the $\alpha$-Perron vector of $G$.

Let $\mathcal{G}$ be a set of graphs. For the work on extremal spectral problems, one of the most important problems is to find the upper or lower bounds for some spectral parameter (index, $Q$-index or $\alpha$-index, etc.) in $\mathcal{G}$ and characterize the extremal graphs.
There are two classic problems related to this problem. One is the Brualdi-Soheild problem \cite{BS}: find an upper bound for the indices in $\mathcal{G}$ of order $n$ and characterize the extremal graphs, and the other is the Brualdi-Hoffman problem \cite{BH}: find an upper bound for the indices in $\mathcal{G}$ of size $m$ and characterize the extremal graphs. For related researches, one may refer to \cite{B,CR,R,Stanic,Stanley,Ste}.

It is interesting to consider the above two problems under the restrictions of other parameters or special classes of graphs. A graph is said to be $H$-free if it does not contain a subgraph isomorphic to $H$. Berman and Zhang \cite{BZ} characterized the graphs with maximum index among all connected graphs with order $n$ and cut vertices $k$. Liu, Lu and Tian \cite{LLT} determined the graphs with the maximum index among all the unicyclic graphs with order $n$ and diameter $d$. Zhai, Lin and Shu \cite{ZLS} characterized the graphs with the maximum index among the $K_{2, r+1}$-free (resp. $\{C_3^+,C_4^+\}$-free) graphs with size $m$. For the $Q$-index and $\alpha$-index counterparts of the above problems, many researchers also have some corresponding results. Zhai, Xue and Lou \cite{ZXL} determined the graph with the maximum $Q$-index among all graphs with size $m$ and clique number $\omega$ (resp.\;chromatic number $\chi$). Lin, Huang and Xue \cite{LHX} characterized the graph with the maximum $\alpha$-index among all connected graphs with order $n$ and cut vertices $k$. Guo and Zhang \cite{GZ1} determined the graphs with the maximum $\alpha$-index for $\alpha \in [\frac{1}{2},1)$ among the $C_4$-free (resp. Halin) graphs with order $n$. For more results, one can refer to \cite{FW,GS,LBW,TCC,ZL}. In recent years, the relationship between the spectral parameter and forbidden subgraphs has been a hot research topic. We refer the interested reader to the surveys \cite{CZ,LLF,N2011} for more results.

A graph is $k$-connected (resp. $k$-edge-connected) if removing fewer than $k$ vertices (resp. edges) always leaves the remaining graph connected, and is minimally $k$-(edge)-connected if it is $k$-connected (resp. $k$-edge-connected) and deleting any arbitrary chosen edge always leaves a graph which is not $k$-connected (resp. $k$-edge-connected). In recent works, some researchers restrict $\mathcal{G}$ to (minimally) $k$-(edge)-connected graphs of order $n$ or size $m$. A graph is minimally 1-(edge)-connected if and only if it is a tree. It is natural to ask which graphs have the maximal indices among all minimally $k$-(edge)-connected graphs for $k\geq 2$. Fan, Goryainov and Lin \cite{FGL} asked the following question for $k\geq 2$.
\noindent\begin{problem}
What is the maximum ($Q$-)index and what are the corresponding extremal graphs among minimally k-(edge)-connected graph for $k\geq 2$?
\end{problem}
Chen and Guo \cite{CG} and Lou, Min and Huang \cite{LMH} characterized the extremal graphs with the maximum index among all minimally 2-(edge)-connected graphs with given order or size, respectively. Fan, Goryainov and Lin \cite{FGL} determined the extremal graphs with the maximum $Q$-index among all minimally 2-(edge)-connected graphs with given order, meanwhile, they characterized the extremal graphs with the maximum ($Q$-)index among all minimally 3-connected graphs with given order. Guo and Zhang \cite{GZ2,ZG} characterized the extremal graphs with the maximum $Q$-index among all (minimally) 2-connected graphs with given size. Analogously, we ask the following question with respect to $\alpha$-index.
\noindent\begin{problem}
What is the maximum $\alpha$-index and what are the corresponding extremal graphs among minimally k-(edge)-connected graph for $k\geq 2$?
\end{problem}
\begin{figure}\label{fig1}
  \centering
  \includegraphics[width=10.5cm]{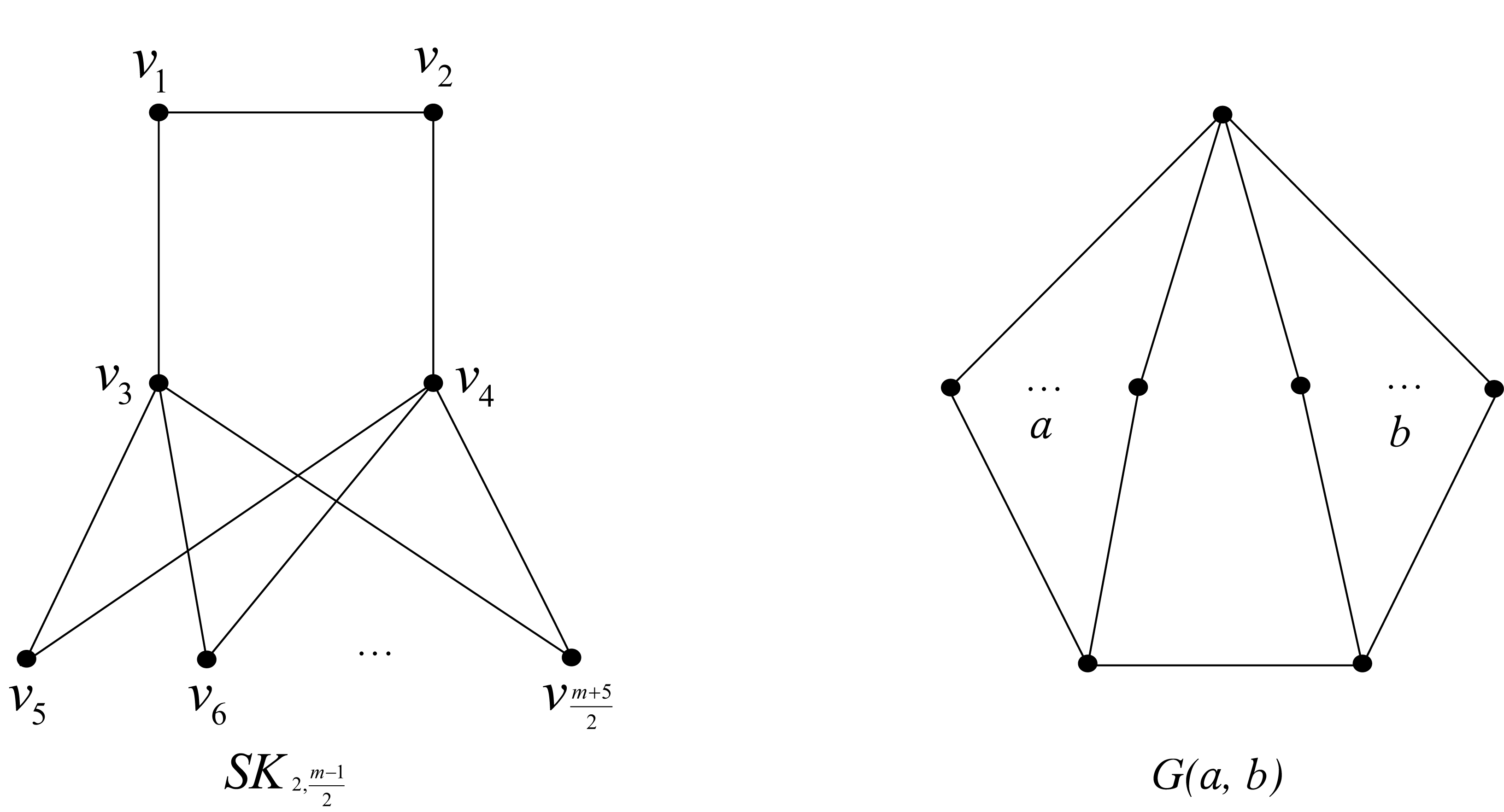}\\
  \small{\textbf{Fig.\;1.} The graphs $SK_{2,\frac{m-1}{2}}$ and $G(a,b)$}
\end{figure}
In this paper, we characterize the extremal graphs with the minimum $\alpha$-index for $\alpha \in [\frac{1}{2},1)$ among all minimally 2-connected graphs with given order or size, respectively.

\noindent\begin{theorem}\label{th1.3}
Let $G$ be a minimally 2-connected graph with order $n\geq5$. If $\alpha \in[\frac{1}{2},1)$, then $\rho_{\alpha}(G)\leq\rho_{\alpha}(K_{2,n-2})$,
with equality if and only if $G\cong K_{2,n-2}$.
\end{theorem}
\noindent\begin{theorem}\label{th1.4}
Let $G$ be a minimally 2-connected graph with size $m$ and $\alpha \in [\frac{1}{2},1)$.
\begin{enumerate}[{\rm(i)}]
\item If $m \geq 6$ is an even number, then $\rho_{\alpha}(G)\leq\rho_{\alpha}(K_{2,\frac{m}{2}})$, with equality if and only if $G\cong K_{2,\frac{m}{2}}$.
\item If $m \geq 9$ is an odd number, then $\rho_{\alpha}(G)\leq \rho_{\alpha}(SK_{2, \frac{m-1}{2}})$, where $\rho_{\alpha}(SK_{2, \frac{m-1}{2}})$ is the largest root of\enspace
$x^3-(\frac{m+5}{2}\alpha+1)x^2+(\frac{m+5}{2}\alpha^2+\frac{5(m-1)}{2}\alpha+2-m)x
 -2m\alpha^2-(m-5)\alpha+m-3=0$,
with equality if and only if $G\cong SK_{2, \frac{m-1}{2}}$.
\end{enumerate}
\end{theorem}
The rest of this paper is organized as follows. In Section 2, we recall some notions and lemmas that will be used later, and prove some new lemmas. In Sections 3 and 4, we give the proof of Theorems \ref{th1.3} and \ref{th1.4}, respectively.

\section{Preliminaries}
\label{sec:Preliminaries}
\noindent In this section, we introduce some preliminary results that are used in the proof of our main results.

\noindent\begin{lemma}\label{lem1} {\normalfont(\cite{N2017})}
If $G$ is a graph with no isolated vertices, then
$$\rho_{\alpha}(G) \leq \max _{u \in V(G)}\left\{\alpha d(u)+\frac{1-\alpha}{d(u)} \sum_{u v \in E(G)} d(v)\right\}.$$
If $\alpha \in (\frac{1}{2},1)$ and $G$ is connected, with equality if and only if $G$ is regular.
\end{lemma}

\noindent\begin{lemma}\label{lem2} {\normalfont(\cite{N2017})}
Let G be a graph with $\Delta(G)=\Delta$. If $\alpha \in[0,\frac{1}{2}]$, then
$$\rho_{\alpha}(G) \geq \alpha(\Delta+1).$$
If $\alpha \in[\frac{1}{2},1)$, then
$$\rho_{\alpha}(G) \geq \alpha \Delta+\frac{(1-\alpha)^{2}}{\alpha}.$$
\end{lemma}

\noindent\begin{lemma}\label{lem3} {\normalfont(\cite{BM})}
If $G$ is a minimally 2-(edge)-connected graph, then $\delta(G)=2 .$
\end{lemma}

\noindent\begin{lemma}\label{lem4} {\normalfont(\cite{D})}
A minimally 2-connected graph with more than three vertices contains no triangles.
\end{lemma}

\noindent\begin{lemma}\label{lem5} {\normalfont(\cite{D})}
A minimally 2-connected graph with $n\geq 4$ has at most $2n-4$ edges, with equality if and only if $G\cong K_{2,n-2}$.
\end{lemma}

\noindent\begin{lemma}\label{lem6} {\normalfont(\cite{P})}
A 2-connected graph $G$ is minimally 2-connected if and only if no cycle of $G$ has a chord.
\end{lemma}

\noindent\begin{lemma}\label{lem7} {\normalfont(\cite{NR,XLL})}
Let $G$ be a connected graph with $\alpha \in[0,1)$. For $u, v \in V(G)$, suppose $N \subseteq N(v) \backslash(N(u) \cup\{u\})$. Let $G^{\prime}=G-\{v w: w \in N\}+\{u w: w \in N\}$. Let $X=\left(x_1, x_2, \ldots, x_n\right)^{T}$ be the $\alpha$-Perron vector of $G$ corresponding to $\rho_\alpha(G)$. If $N \neq \emptyset$ and $x_u \geq x_v$, then $\rho_\alpha\left(G^{\prime}\right)>\rho_\alpha(G)$.
\end{lemma}
We say that $u$ and $v$ are equivalent in $G$ if there exists an automorphism $p: G\rightarrow G$ such that $p(u)=v$.
\noindent\begin{lemma}\label{lem8} {\normalfont(\cite{N2017})}
Let $G$ be a connected graph of order $n$, and let $u$ and $v$ be equivalent vertices in G. If $X=\left(x_{1}, x_{2}, \ldots, x_{n}\right)^{T}$ is an eigenvector to $\rho_{\alpha}(G)$, then $x_{u}=x_{v}$.
\end{lemma}

\noindent\begin{lemma}\label{lem9}
{\normalfont(\cite{N2017})}
Let $a \geq b \geq 1$. If $\alpha \in[0,1]$, then the largest eigenvalue of $A_\alpha\left(K_{a, b}\right)$ is
$$
\rho_\alpha(K_{a,b})=\frac{1}{2}(\alpha(a+b)+\sqrt{\alpha^2(a+b)^2+4 a b(1-2 \alpha)}).
$$
\end{lemma}

\noindent\begin{lemma}\label{lem10}
Let $m$ be an odd integer. Then $\rho_ \alpha{ (SK_{2,\frac{m-1}{2}})}$ is the largest root of the following equation:
$$
\begin{aligned}
x^3-(\frac{m+5}{2}\alpha+1)x^2+(\frac{m+5}{2}\alpha^2+\frac{5(m-1)}{2}\alpha+2-m)x
-2m\alpha^2-\left(m-5\right )\alpha+m-3=0.
\end{aligned}
$$
\end{lemma}
\noindent\textbf{Proof.}
Let $V(SK_{2,\frac{m-1}{2}})=\{v_1,v_2,...,v_\frac{m+5}{2}\}$ (see Fig.\;1) and $X=(x_1,x_2,...,x_\frac{m+5}{2})^{T}$ be the $\alpha$-Perron vector of $SK_{2,\frac{m-1}{2}}$, where $x_i$ denotes the coordinate corresponding to $v_i$ for $1\leq i\leq\frac{m+5}{2}$. By Lemma \ref{lem8}, we have
$$x_1=x_2,\ x_3=x_4,\ x_5=...=x_{\frac{m+5}{2}}.$$
 Let $\rho_\alpha=\rho_\alpha{(SK_{2,\frac{m-1}{2}})}$. Since $A_\alpha(SK_{2,\frac{m-1}{2}})X=\rho_\alpha X$, then
$$
\left\{
\begin{aligned}
\left(\rho_\alpha-\left (\alpha+1\right )\right) x_{1} &=(1-\alpha) x_{3}, \\
(\rho_\alpha-\frac{(m-1)\alpha}{2})\, x_{3} &=(1-\alpha)x_{1}+\frac{(m-3)(1-\alpha)}{2} x_{5}, \\
\left(\rho_\alpha-2\alpha\right ) x_{5} &=2(1-\alpha) x_{3}.
\end{aligned}
\right.
$$
Since $X=(x_{1}, x_{2}, \ldots, x_{\frac{m+5}{2}})^{T}$ is an eigenvector corresponding to $\rho_\alpha$, it follows that $X \neq 0$. This implies that
$$
\left|\begin{array}{ccc}
\rho_\alpha-(\alpha+1) & -(1-\alpha) & 0 \\
-(1-\alpha) & \rho_\alpha-\frac{(m-1)\alpha}{2} & -\frac{(m-3)(1-\alpha)}{2} \\
0 & -2(1-\alpha) & \rho_\alpha-2\alpha
\end{array}\right|=0.
$$
Hence $\rho_\alpha$ is the largest root of the following equation
$$
\left|\begin{array}{ccc}
x-(\alpha+1) & -(1-\alpha) & 0 \\
-(1-\alpha) & x-\frac{(m-1)\alpha}{2} & -\frac{(m-3)(1-\alpha)}{2} \\
0 & -2(1-\alpha) & x-2\alpha
\end{array}\right|=0.
$$
By computation, we conclude that $\rho_\alpha$ is the largest root of the following equation:
$$
\begin{aligned}
x^3-(\frac{m+5}{2}\alpha+1)x^2+(\frac{m+5}{2}\alpha^2+\frac{5(m-1)}{2}\alpha+2-m)x-2m\alpha^2-(m-5)\alpha+m-3=0.
\end{aligned}
$$

This completes the proof.\ $\qedsymbol$
\noindent\begin{lemma}\label{lem11}
Let $m\geq 9$ and $\alpha \in [\frac{1}{2}, 1)$. Then
$$f(\alpha,m)>0 \ and\  g(\alpha,m)<0,$$
$
\begin{aligned}where \enspace
f(\alpha,m)=&\, 2(5\alpha^3-6\alpha^2+2\alpha)m^5-2(123\alpha^3-156\alpha^2+60\alpha-4)m^4
+4(517\alpha^3-650\alpha^2\\&+254\alpha-20)m^3-4(1783\alpha^3-2008\alpha^2+664\alpha-16)m^2
+2(5377\alpha^3-5014\alpha^2\\&+1202\alpha+136)m-2(2951\alpha^3-2052\alpha^2+228\alpha+196)
\end{aligned}
$\\
$
\begin{aligned} and \enspace
g(\alpha,m)=&\, (2\alpha^4+6\alpha^3-9\alpha^2+1)m^3+(8\alpha^4+9\alpha^3-34\alpha^2+\alpha)m^2
-2(35\alpha^4+154\alpha^3\\&-135\alpha^2+20\alpha+14)m-4(75\alpha^4-79\alpha^3+137\alpha^2-21\alpha-16).
\end{aligned}
$
\end{lemma}
\noindent\textbf{Proof.}
The first, second and third order partial derivatives of function $f(\alpha, m)$ with respect to $m$ are as follows:
$$
\begin{aligned}
f_{m}(\alpha, m)=&\, 10(5\alpha^3-6\alpha^2+2\alpha)m^4-8(123\alpha^3-156\alpha^2+60\alpha-4)m^3
+12(517\alpha^3-650\alpha^2\\&+254\alpha-20)m^2-8(1783\alpha^3-2008\alpha^2+664\alpha-16)m
+2(5377\alpha^3-5014\alpha^2\\&+1202\alpha+136), \\
f_{mm}(\alpha, m)=&\, 40(5\alpha^3-6\alpha^2+2\alpha)m^3-24(123\alpha^3-156\alpha^2+60\alpha-4)m^2
+24(517\alpha^3-650\alpha^2\\&+254\alpha-20)m-8(1783\alpha^3-2008\alpha^2+664\alpha-16), \\
f_{mmm}(\alpha,m)=&\, 120(5\alpha^3-6\alpha^2+2\alpha)m^2-48(123\alpha^3-156\alpha^2+60\alpha-4)m
+24(517\alpha^3-650\alpha^2\\&+254\alpha-20). \\
\end{aligned}
$$
Since $\alpha \in [\frac{1}{2}, 1)$, then $5\alpha^{3}-6 \alpha^{2}+2\alpha>0$. Let $m_0$ be the minimum point of $f_{mmm}(\alpha,m)$. Then we have
$$
m_{0}=\frac{48(123 \alpha^{3}-156\alpha^{2}+60\alpha-4)}{240\left(5 \alpha^{3}+6 \alpha^{2}-2\alpha\right)}=\frac{123\alpha^3-156 \alpha^{2}+60\alpha-4}{5(5\alpha^3-6 \alpha^2+2 \alpha)}<9.
$$
It follows that $f_{mmm}(\alpha, m)$ is increasing for $m \geq 9$ and
$$
f_{mmm}(\alpha, m) \geq f_{mmm}(\alpha, 9)=96(82\alpha^{3}-68\alpha^2-4\alpha+13).
$$
It is easy to verify that $82\alpha^{3}-68\alpha^2-4\alpha+13>0$ for $\alpha \in [\frac{1}{2}, 1)$, that is, $f_{mmm}(\alpha, m)>0$. It follows that $f_{mm}(\alpha, m)$ is increasing for $m \geq 9$ and
$$
f_{mm}(\alpha, m) \geq f_{mm}(\alpha,9)=64(64\alpha^3+62\alpha^2-137\alpha+56).
$$
Note that $64\alpha^3+62\alpha^2-137\alpha+56>0$ for $\alpha \in [\frac{1}{2}, 1)$, that is, $f_{mm}(\alpha, m)>0$. It follows that $f_{m}(\alpha, m)$ is increasing for $m \geq 9$ and
$$
f_{m}(\alpha, m) \geq f_{m}(\alpha,9)=-32(137\alpha^3-590\alpha^2+538\alpha-166).
$$
It is obvious that $137\alpha^3-590\alpha^2+538\alpha-166<0$ for $\alpha \in [\frac{1}{2}, 1)$, that is, $f_m(\alpha, m)>0$. It follows that $f(\alpha, m)$ is increasing for $m \geq 9$ and
$$
f(\alpha, m) \geq f(\alpha, 9)=-64(43\alpha^3-117\alpha^2+69\alpha-22)>0.
$$

Similarly, we have $$
\begin{aligned}
g(\alpha,m)=&\, (2\alpha^4+6\alpha^3-9\alpha^2+1)m^3+(8\alpha^4+9\alpha^3-34\alpha^2+\alpha)m^2
-2(35\alpha^4+154\alpha^3-135\alpha^2\\&+20\alpha+14)m-4(75\alpha^4-79\alpha^3+137\alpha^2-21\alpha-16)<0.
\end{aligned}
$$ for $m\geq9$ and $\alpha \in [\frac{1}{2}, 1)$.

This completes the proof.\ $\qedsymbol$


\section{Proof of Theorem \ref{th1.3}}
\label{sec:ch-sufficient}
\noindent In this section, we give the proof of Theorem \ref{th1.3}.\\ \hspace*{\fill} \\
\textbf{Proof.}
Let $G$ be a minimally 2-connected graph with order $n\geq5$, then $\Delta(G)\leq n-2$ by Lemmas \ref{lem3} and \ref{lem4}. Notice that $K_{2,n-2}$ is a minimally 2-connected graph. By Lemma \ref{lem9}, we have
$$
\begin{aligned}
\rho_{\alpha}(K_{2,n-2})=\frac{1}{2}(\alpha n+\sqrt{\alpha^2 n^2+8(n-2)(1-2\alpha)}).
\end{aligned}
$$
When $\Delta(G)=n-2$, it is easy to see that $G\cong K_{2,n-2}$ by Lemmas \ref{lem3} and \ref{lem4}. By Lemma \ref{lem2}, we have
\begin{equation}{\label{1}}
\rho_{\alpha}(K_{2,n-2})\geq \alpha\Delta(K_{2,n-2})+\frac{(1-\alpha)^{2}}{\alpha}
=\alpha(n-2)+\frac{(1-\alpha)^{2}}{\alpha}
\end{equation}
for $\alpha\in[\frac{1}{2},1)$. Thus we assume that $\Delta(G)\leq n-3$.\\
Let $w$ be a vertex of $G$ such that
$$\alpha d(w)+\frac{1-\alpha}{d(w)} \sum_{wv \in E(G)}d(v)=\max_{u \in V(G)}\left\{\alpha d(u)+\frac{1-\alpha}{d(u)} \sum_{uv \in E(G)} d(v)\right\}.$$
By Lemma \ref{lem1}, we have
\begin{equation}{\label{2}}
\rho_{\alpha}(G)\leq \alpha d(w)+\frac{1-\alpha}{d(w)} \sum_{wv \in E(G)}d(v).
\end{equation}
By Lemma \ref{lem3}, we have $2\leq d(w) \leq \Delta(G)\leq n-3$. By Lemma \ref{lem4}, we know that $N(w)$ is an independent set, that is, $e(N(w))=0$. Then
\begin{equation}{\label{3}}
\begin{aligned}
 \sum_{wv \in E(G)} d(v)=2 e(N(w))+e(N(w), V(G) \backslash N(w))=e(N(w), V(G) \backslash N(w)).
 \end{aligned}
\end{equation}

Next we prove that $\rho_\alpha(G)\leq \rho_\alpha(K_{2,n-2})$ for $2\leq d(w)\leq n-3$. We consider the following two cases.

\textbf{Case 1.} $d(w)=2$.

If $e(V(G)\backslash N[w])=0$, then $G\cong K_{2,n-2}$ by Lemmas \ref{lem3} and \ref{lem4}.

If $e(V(G)\backslash N[w])\neq0$. We assume that there exists an edge $v_1v_2\in E(G[V(G)\backslash N[w]])$. By Lemma \ref{lem4}, we have $N_N(w)(v_1)\cap N_N(w)(v_1)=\varnothing$. Let $B=V(G) \backslash (N[w]\cup \{v_1,v_2\})$. Then
$$e(N(w), V(G) \backslash N[w])\leq d_{N(w)}(v_1)+d_{N(w)}(v_2)+ e(N(w),B)= d(w)+d(w)|B|=2n-8.$$
and so $\sum_{wv \in E(G)} d(v)\leq 2n-6$. Combining this with (\ref{2}), we have
$$
\rho_\alpha(G) \leq 2 \alpha+\frac{1-\alpha}{2} \sum_{w v \in E(G)} d(v) \leq 2 \alpha+(1-\alpha)(n-3).
$$
Noting that
$$
\alpha(n-2)+\frac{(1-\alpha)^2}{\alpha}-(2 \alpha+(1-\alpha)(n-3))=\frac{(2\alpha^2-\alpha)n-6 \alpha^2+ \alpha+1}{ \alpha}\geq0
$$
for $n \geq 5$ and $\alpha \in[\frac{1}{2},1)$, we have
$\rho_\alpha(G)\leq\alpha(n-2)+\frac{(1-\alpha)^2}{\alpha}.$
Combining this with (\ref{1}), we have
$$
\rho_\alpha(G)\leq\alpha(n-2)+\frac{(1-\alpha)^2}{\alpha}\leq\rho_\alpha\left(K_{2, n-2}\right)
$$
for $n \geq 5$ and $\alpha \in[\frac{1}{2},1)$.

\textbf{Case 2.} $3\leq d(w) \leq n-3$.

In order to prove $\rho_\alpha(G)\leq\rho_\alpha(K_{2,n-2})$, it is enough to prove
$$
\rho_\alpha(G)\leq \frac{1}{2}(\alpha n+\sqrt{\alpha^2 n^2+8(n-2)(1-2\alpha)}),
$$
that is, to prove
$\rho_\alpha(G)^2-\alpha n \rho_\alpha(G) +2(2\alpha-1)(n-2)\leq 0.$
For convenience, we denote $A_\alpha(G)=A_\alpha$, $A(G)=A$ and $D(G)=D$. Let
$$
B=(b_{ij})_{n\times n}=A_\alpha^2-\alpha n A_\alpha +2(2\alpha-1)(n-2)I_n,
$$
where $I_n$ is the $n\times n$ unit matrix. Let $c_u(B)$ be the sum of all elements in the $u$-th column of $B$. Then we have the following claim.
\par{\textbf{Claim 2.1.} $c_u(B)\leq 0$ for $n\geq 5$ and $\alpha \in[\frac{1}{2},1)$.
\par{\textbf{Proof.}
Since $A_\alpha=\alpha D+(1-\alpha)A$, then
$$
\begin{aligned}
B=&(\alpha D+(1-\alpha) A)^2-\alpha n(\alpha D+(1-\alpha) A)+2(2\alpha-1)(n-2) I_n \\
=& \alpha^2 D^2+(1-\alpha)^2 A^2+\alpha(1-\alpha) D A+\alpha(1-\alpha) A D-\alpha^2 n D
-\left(\alpha n-\alpha^2 n\right) A\\&+2(2\alpha-1)(n-2)I_n.
\end{aligned}
$$

It is easy to see that $c_u(A)=c_u(D)=d(u)$, $c_u(A^2)=c_u(DA)=\sum_{u v \in E(G)} d(v)$ and $c_u(AD)=d^2(u)$. Combining (\ref{3}) with Lemma \ref{lem5}, we have
$$
\begin{aligned}
\sum_{w v \in E(G)} d(v)=e(N(w), V(G) \backslash N(w))\leq |E(G)|\leq 2n-4,
\end{aligned}
$$
with equality if and only if $G\cong K_{2,n-2}$. It follows that
$$
\begin{aligned}
c_u(B)=& \alpha^2 d^2(u)+(1-\alpha)^2 \sum_{uv \in E(G)} d(v)+\alpha(1-\alpha) \sum_{uv \in E(G)} d(v)+\alpha(1-\alpha)d^2(u)-\alpha^2 n d(u)
\\&-\left(\alpha n-\alpha^2 n\right) d(u)+2(2\alpha-1)(n-2) \\
=&\alpha^2 d^2(u)+(1-\alpha) \sum_{uv \in E(G)} d(v)-\alpha n d(u)+2(2\alpha-1)(n-2)\\
\leq & \alpha^2 d^2(u)+(1-\alpha) (2n-4)-\alpha n d(u)+2(2\alpha-1)(n-2) \\
=& \alpha\left(d^2(u)-n d(u)+2 n-4\right) \\
\leq & \max \left\{\alpha(9-3n+2 n-4), \alpha\left((n-3)^2-n(n-3)+2 n-4\right)\right\} \\
=& \alpha(-n+5)\leq 0
\end{aligned}
$$
for $n\geq 5$ and $\alpha \in[\frac{1}{2},1)$, with equality if and only if $G\cong K_{2,n-2}$.

This completes the proof of the claim. $\qedsymbol$}
}
\par{Let $X=\left(x_1, x_2, \ldots, x_n\right)^{T}$ be the $\alpha$-Perron vector of $G$ corresponding to $\rho_\alpha(G)$ satisfying $\sum_{i=1}^n x_i=1$. Then
$$
B X=\left(\rho_\alpha(G)^2-\alpha n \rho_\alpha(G) +2(2\alpha-1)(n-2)\right)X.
$$
Hence we have
$$
\begin{aligned}
&\rho_\alpha(G)^2-\alpha n \rho_\alpha(G) +2(2\alpha-1)(n-2)\\
&=\sum_{i=1}^n\left(\rho_\alpha(G)^2-\alpha n \rho_\alpha(G) +2(2\alpha-1)(n-2)\right) x_i \\
&=\sum_{i=1}^n(B X)_i=\sum_{i=1}^n(\sum_{j=1}^n b_{i j} x_j) =\sum_{j=1}^n(\sum_{i=1}^n b_{i j}) x_j=\sum_{j=1}^n c_j(B) x_j\leq0 .
\end{aligned}
$$}
\par{Combining the above arguments, we have $\rho_\alpha(G)\leq \rho_\alpha(K_{2,n-2})$ for $n\geq 5$ and $\alpha \in[\frac{1}{2},1)$, with equality if and only if $G\cong K_{2,n-2}$.

These complete the proof. $\qedsymbol$}

\section{Proof of Theorem \ref{th1.4}}
\noindent In this section, we give the proof of Theorem \ref{th1.4}.\\ \hspace*{\fill} \\
\textbf{Proof.}
Let $G$ be a minimally 2-connected graph with size $m$. For any $v \in V(G)$, it is easy to see that $G-v$ is connected and $|E(G-v)| = m-d(v)$, then we have
$$
\begin{aligned}
d(v)< |V(G-v)|\leq m-d(v)+1,
\end{aligned}
$$
with equality in the right inequality if and only if $G-v$ is a tree. It follows that $d(v)< \frac{m+1}{2}$. Combining this with Lemma \ref{lem3}, we have $ 2\leq d(v)< \frac{m+1}{2}$.\\
Let $w$ be a vertex of $G$ such that
$$\alpha d(w)+\frac{1-\alpha}{d(w)} \sum_{wv \in E(G)}d(v)=\max_{u \in V(G)}\left\{\alpha d(u)+\frac{1-\alpha}{d(u)} \sum_{uv \in E(G)} d(v)\right\}.$$
By Lemma \ref{lem1}, we have
\begin{equation}{\label{4}}
\rho_{\alpha}(G)\leq \alpha d(w)+\frac{1-\alpha}{d(w)} \sum_{wv \in E(G)}d(v).
\end{equation}
By Lemma \ref{lem4}, we know that $N(w)$ is an independent set, that is, $e(N(w))=0$. Then
\begin{equation}{\label{5}}
\begin{aligned}
\sum_{wv \in E(G)}d(v)=2e(N(w))+e(N(w), V(G) \backslash N(w))=e(N(w), V(G) \backslash N(w)).
\end{aligned}
\end{equation}
Combining (\ref{4}) and (\ref{5}), we have
\begin{equation}{\label{6}}
\begin{aligned}
\rho_{\alpha}(G)\leq \alpha d(w)+\frac{1-\alpha}{d(w)}e(N(w), V(G) \backslash N(w)).
\end{aligned}
\end{equation}

\textbf{(\romannumeral1)} Let $m\geq 6$ be an even number. Then we have $2\leq d(w)\leq \frac{m}{2}$. Notice that $K_{2,\frac{m}{2}}$ is a minimally 2-connected graph. By Lemma \ref{lem9}, we have
$$
\begin{aligned}
\rho_\alpha(K_{2,\frac{m}{2}})
=\frac{1}{4}((m+4)\alpha+\sqrt{(m+4)^2\alpha^2+16m(1-2\alpha)}).
\end{aligned}
$$
If $d(w)=\frac{m}{2}$, then $d(v)=2$ for any $v \in N(w)$ and $e(N(w),V(G)\backslash N[w])=\frac{m}{2}$ by Lemma \ref{lem3}. Since $G-w$ is connected, then $G\cong K_{2, \frac{m}{2}}$.\\
Next we prove that $\rho_\alpha(G)\leq \rho_\alpha(K_{2, \frac{m}{2}})$ for $2\leq d(w)\leq\frac{m-2}{2}$. We consider the following two cases.

\textbf{Case 1.} $d(w)=2$.

If $e(V(G)\backslash N[w])=0$, then $G\cong K_{2,\frac{m}{2}}$ by Lemmas \ref{lem3} and \ref{lem4}.

If $e(V(G)\backslash N[w])\neq0$, then
$e(N(w), V(G) \backslash N(w))\leq m-1$ by Lemma \ref{lem4}. Combining this with (\ref{6}), we have
$$
\rho_\alpha(G) \leq 2 \alpha+\frac{1-\alpha}{2}(m-1)\leq \frac{m+3}{4}<\frac{m+4}{4}
$$
for $m \geq 6$ and $\alpha \in[\frac{1}{2}, 1)$.
In order to prove $\rho_\alpha(G)<\rho_\alpha(K_{2,\frac{m}{2}})$, it is enough to prove
$$
\frac{m+4}{4}\leq \frac{1}{4}((m+4)\alpha+\sqrt{(m+4)^2\alpha^2+16m(1-2\alpha)}),
$$
that is, to prove $(1-2\alpha)m^2-8(1-2\alpha)m+16(1-2\alpha)\leq 0$
for $m \geq 6$ and $\alpha \in[\frac{1}{2}, 1)$. It is easy to check that this is true. Hence we have
$\rho_\alpha(G)<\rho_\alpha(K_{2, \frac{m}{2}})$
for $m \geq 6$ and $\alpha \in[\frac{1}{2},1)$.

\textbf{Case 2.} $3\leq d(w) \leq \frac{m-2}{2}$.

In order to prove $\rho_\alpha(G)\leq\rho_\alpha(K_{2,\frac{m}{2}})$, it is enough to prove
$$
\rho_\alpha(G)\leq \frac{1}{4}((m+4)\alpha+\sqrt{(m+4)^2\alpha^2+16m(1-2\alpha)}),
$$
that is, to prove
$2\rho_\alpha(G)^2-(m+4)\alpha \rho_\alpha(G) +2(2\alpha-1)m<0.$
For convenience, we denote $A_\alpha(G)=A_\alpha$, $A(G)=A$ and $D(G)=D$. Let
$$
B=(b_{ij})_{n\times n}=2A_\alpha^2-(m+4)\alpha A_\alpha +2(2\alpha-1)mI_n,
$$
where $I_n$ is the $n\times n$ unit matrix. Let $c_u(B)$ be the sum of all elements in the $u$-th column of matrix $B$. Then we have the following claim.
\par{\textbf{Claim 2.1.} $c_u(B)\leq 0$ for $m\geq 6$ and $\alpha \in[\frac{1}{2},1)$.
\par{\textbf{Proof.}
Since $A_\alpha=\alpha D+(1-\alpha)A$, then
$$
\begin{aligned}
B=&2(\alpha D+(1-\alpha) A)^2-(m+4)\alpha (\alpha D+(1-\alpha) A)+2(2\alpha-1)m I_n \\
=& 2\alpha^2 D^2+2(1-\alpha)^2 A^2+2\alpha(1-\alpha) D A+2\alpha(1-\alpha) A D-(m+4)\alpha^2 D
\\&-(m+4)\alpha(1-\alpha) A+2(2\alpha-1)mI_n.
\end{aligned}
$$
It is easy to see that $c_u(A)=c_u(D)=d(u)$, $c_u(A^2)=c_u(DA)=\sum_{u v \in E(G)} d(v)$ and $c_u(AD)=d^2(u)$. Since $e(N(w), V(G) \backslash N(w))\leq |E(G)|=m$, then $\sum_{wv \in E(G)}d(v) \leq m$ by (\ref{5}). It follows that
$$
\begin{aligned}
c_u(B)=& 2\alpha^2 d^2(u)+2(1-\alpha)^2 \sum_{wv \in E(G)} d(v)+2\alpha(1-\alpha) \sum_{wv \in E(G)} d(v)+2\alpha(1-\alpha)d^2(u)\\&-(m+4)\alpha^2 d(u)
-(m+4)\alpha(1-\alpha) d(u)+2(2\alpha-1)m \\
=&2\alpha^2 d^2(u)+2(1-\alpha) \sum_{wv \in E(G)} d(v)-(m+4)\alpha d(u)+2(2\alpha-1)m\\
\leq & 2\alpha d^2(u)-(m+4)\alpha d(u)+2m\alpha
\end{aligned}
$$
$$
\begin{aligned}
=& \alpha(2d^2(u)-(m+4) d(u)+2m) \\
\leq & \max \left\{\alpha(18-3(m+4)+2m), \alpha\left(\frac{(m-2)^2}{2}-\frac{(m+4)(m-2)}{2}+2m\right)\right\} \\
=& \alpha(-m+6)\leq 0
\end{aligned}
$$
for $m\geq 6$ and $\alpha \in[\frac{1}{2},1)$, with equality if and only if $G\cong K_{2,\frac{m}{2}}$.

This completes the proof of the claim. $\qedsymbol$}
}
\par{Let $X=\left(x_1, x_2, \ldots, x_n\right)^{T}$ be the $\alpha$-Perron vector of $G$ corresponding to $\rho_\alpha(G)$ satisfying $\sum_{i=1}^n x_i=1$. Then
$$
BX=\left(2\rho_\alpha(G)^2-(m+4)\alpha \rho_\alpha(G) +2(2\alpha-1)m\right)X.
$$
Hence we have
$$
\begin{aligned}
&2\rho_\alpha(G)^2-(m+4)\alpha \rho_\alpha(G) +2(2\alpha-1)m\\
&=\sum_{i=1}^n\left(2\rho_\alpha(G)^2-(m+4)\alpha \rho_\alpha(G) +2(2\alpha-1)m\right) x_i \\
&=\sum_{i=1}^n(BX)_i=\sum_{i=1}^n(\sum_{j=1}^n b_{i j} x_j) =\sum_{j=1}^n(\sum_{i=1}^n b_{i j}) x_j=\sum_{j=1}^n c_j(B) x_j\leq0.
\end{aligned}
$$}
\par{Combining the above arguments, we have $\rho_{\alpha}(G) \leq \rho_{\alpha}(K_{2,\frac{m}{2}})$
for $m\geq 6$ and $\alpha \in [\frac{1}{2}, 1)$, with equality if and only if $G\cong K_{2,\frac{m}{2}}$. We complete the proof of (i).}

\textbf{(\romannumeral2)} Let $m\geq 9$ be an odd number. Then we have $2\leq d(w)\leq \frac{m-1}{2}$. Notice that $SK_{2,\frac{m}{2}}$ is a minimally 2-connected graph. Next we complete the proof with three facts.

\textbf{Fact 1.} $e(N(w), V(G)\backslash N(w))\leq m-1$.

\textbf{Proof.}
Since $e(N(w), V(G) \backslash N(w))\leq |E(G)|=m$, then $\sum_{wv \in E(G)}d(v) \leq m$ by (\ref{5}). For a contradiction, we suppose $e(N(w), V(G)\backslash N(w))= m$, then $e(V(G)\backslash N[w])= 0$. We consider the following two cases.

\textbf{Case 1.1.} $d(w)=2$.

In this case, we have $d(v)=2$ for any $v \in V(G)\backslash N(w)$ by Lemma \ref{lem3}. This implies that $G$ is a complete bipartite graph $K_{2,b}$ and $m(K_{2,b})=2b$, which contradicts the fact that $m$ is odd.

\textbf{Case 1.2.} $3 \leq d(w) \leq \frac{m-1}{2}$.

Let $v_1, v_2 \in V(G)\backslash N[w]$ be any two vertices. If $N_{N(w)}(v_1)\cap N_{N(w)}(v_2) = \varnothing$ , then $d(w_i)=2$ for any $w_i \in N(w)$ by Lemma \ref{lem3}. It follows that $m=e(N(w),V(G)\backslash N(w))=2d(w)\leq m-1$, a contradiction.
If $N_{N(w)}(v_1)\cap N_{N(w)}(v_2) \neq \varnothing$ and $N_{N(w)}(v_1)\neq N_{N(w)}(v_2)$, we assume $w_{12} \in N_{N(w)}(v_1)\cap N_{N(w)}(v_2)$. By Lemma \ref{lem3}, there exists $w_i \in N_{N(w)}(v_i)\backslash w_{12}$ for each $i \in \{1,2\}$. Hence, $G$ contains a cycle $ww_1v_1w_{12}v_2w_2w$ with a chord $ww_{12}$, which contradicts Lemma \ref{lem6}.
If $N_{N(w)}(v_1)= N_{N(w)}(v_2) = N(w)$, then we have $\delta(G)\geq 3$. This contradicts with Lemma \ref{lem3}. Hence $N_{N(w)}(v_1)= N_{N(w)}(v_2) \neq N(w)$. It follows that $G-w$ is disconnected, a contradiction.

Through the above two cases, we know that $e(N(w), V(G)\backslash N(w))\leq m-1$.\ $\qedsymbol$

\textbf{Fact 2.} $\rho_{\alpha}(G) < \rho_{\alpha}{(SK_{2,\frac{m-1}{2}})}$ for $3\leq d(w)\leq \frac{m-3}{2}$.

\textbf{Proof.}
Combining Fact 1 and (\ref{6}), we have
$$
\begin{aligned}
\rho_{\alpha}(G) \leq \alpha d(w)+\frac{1-\alpha}{d(w)} (m-1).
\end{aligned}
$$
Let $q(x)=\alpha x+\frac{1-\alpha}{x}(m-1)$. Since $\alpha \in [\frac{1}{2},1)$, it is easy to see that the function $q(x)$ is convex for $x>0$ and its maximum in any closed internal is attained at one of the ends of this internal. Hence when $3\leq x\leq\frac{m-3}{2}$, we have
$$
\begin{aligned}
\rho_{\alpha}(G) &\leq \max \left\{3 \alpha+\frac{1-\alpha}{3}(m-1),
\frac{m-3}{2} \alpha+\frac{2(1-\alpha)}{m-3}(m-1)\right\} .
\end{aligned}
$$
Noting that
$$
\begin{aligned}
&\frac{m-3}{2}\alpha+\frac{2(1-\alpha)}{m-3}(m-1)-
(3\alpha+\frac{1-\alpha}{3}(m-1)) \\
& = \frac{(5\alpha-2)m^2-(56\alpha-20)m+99\alpha-18}{6(m-3)} > 0
\end{aligned}
$$
for $m\geq 9$ and $\alpha \in [\frac{1}{2},1)$, we have
$$\rho_\alpha(G) \leq \frac{m-3}{2}\alpha+\frac{2(1-\alpha)}{m-3}(m-1)$$
for $m\geq 9$ and $\alpha \in [\frac{1}{2},1)$.
By Lemma \ref{lem10}, we know
$\rho_{\alpha}{(SK_{2,\frac{m-1}{2}})}$ is the largest root of the following equation:
$$
\begin{aligned}
p(x)=&\,x^3-(\frac{m+5}{2}\alpha+1)x^2+(\frac{m+5}{2}\alpha^2+\frac{5(m-1)}{2}\alpha+2-m)x
\\&-2m\alpha^2-(m-5)\alpha+m-3=0.
\end{aligned}
$$
By computation, we have
$$-8(m-3)^3\cdot p(\frac{m-3}{2}\alpha+\frac{2(1-\alpha)}{m-3}(m-1))=f(\alpha,m),$$
$$
\begin{aligned}\textrm{where}  \enspace
f(\alpha,m)=&\,  2(5\alpha^3-6\alpha^2+2\alpha)m^5-2(123\alpha^3-156\alpha^2+60\alpha-4)m^4
+4(517\alpha^3-650\alpha^2\\&+254\alpha-20)m^3-4(1783\alpha^3-2008\alpha^2+664\alpha-16)m^2
+2(5377\alpha^3-5014\alpha^2\\&+1202\alpha+136)m-2(2951\alpha^3-2052\alpha^2+228\alpha+196).
\end{aligned}
$$
By Lemma \ref{lem11}, we have $f(\alpha,m)>0$ for $m\geq 9$ and $\alpha \in [\frac{1}{2}, 1)$. It follows that
$$ p(\frac{m-3}{2}\alpha+\frac{2(1-\alpha)}{m-3}(m-1))<0$$ for $\alpha \in [\frac{1}{2}, 1)$ and $m\geq 9$. This implies that
$$
\begin{aligned}
&\rho_{\alpha}(G) \leq  \frac{m-3}{2}\alpha+\frac{2(1-\alpha)}{m-3}(m-1)<\rho_{\alpha}(SK_{2,\frac{m-1}{2}})
\end{aligned}
$$
for $m\geq 9$ and $\alpha \in [\frac{1}{2},1)$. $\qedsymbol$

\textbf{Fact 3.} $\rho_{\alpha}(G) \leq \rho_{\alpha}{(SK_{2,\frac{m-1}{2}})}$ for $d(w)=2$ or $d(w)= \frac{m-1}{2}$, with equality if and only if $G\cong SK_{2,\frac{m-1}{2}}$.

\textbf{Proof.}
We consider the following two cases.

\par{\textbf{Case 3.1.} $d(w)=2$.}

We consider the following two cases and assume $N(w)=\{w_1, w_2\}$.
\par{\textbf{Subcase 3.1.1.} $e(N(w), V(G)\backslash N(w))= m-1.$
\par{Since $e(N(w), V(G)\backslash N(w))= m-1$, then $e(V(G)\backslash N[w])= 1$. Let $e(V(G)\backslash N[w])={v_1v_2}$. By Lemma \ref{lem4}, we can see that $w_1(resp.\;w_2)$ is adjacent to only one vertex of $v_1$ and $v_2$. Without loss of generality, we assume that $w_1v_1 \in E(G)$ and $w_2v_2 \in E(G)$. By Lemma \ref{lem3}, we have $N_{N(w)}(v)=\{w_1, w_2\}$ for any $v \in V(G)\backslash (N[w]\cup\{v_1,v_2\})$. It follows that $G\cong SK_{2,\frac{m-1}{2}}$.}}
\par{\textbf{Subcase 3.1.2.} $e(N(w), V(G)\backslash N(w))\leq m-2.$
\par{By (\ref{6}), we have
$$\rho_\alpha(G)\leq 2\alpha+\frac{1-\alpha}{2}(m-2).$$
By Lemma \ref{lem10}, we know
$\rho_{\alpha}{(SK_{2,\frac{m-1}{2}})}$ is the largest root of the following equation:
$$
\begin{aligned}
p(x)=&\,x^3-(\frac{m+5}{2}\alpha+1)x^2+(\frac{m+5}{2}\alpha^2+\frac{5(m-1)}{2}\alpha+2-m)x\\
&-2m\alpha^2-(m-5)\alpha+m-3=0.
\end{aligned}
$$
By computation, we have
$$4\cdot p(2\alpha+\frac{1-\alpha}{2}(m-2))=g(\alpha,m),$$
$
\begin{aligned}\textrm{where} \enspace
g(\alpha,m)=&\, (2\alpha^4+6\alpha^3-9\alpha^2+1)m^3+(8\alpha^4+9\alpha^3-34\alpha^2+\alpha)m^2
-2(35\alpha^4+154\alpha^3\\&-135\alpha^2+20\alpha+14)m-4(75\alpha^4-79\alpha^3+137\alpha^2-21\alpha-16).
\end{aligned}
$\\
By Lemma \ref{lem11}, we have $g(\alpha,m)<0$ for $m\geq 9$ and $\alpha \in [\frac{1}{2}, 1)$. It follows that
$$ p(2\alpha+\frac{1-\alpha}{2}(m-2))<0$$ for $m\geq 9$ and $\alpha \in [\frac{1}{2}, 1)$. This implies that
$$
\begin{aligned}
&\rho_{\alpha}(G) \leq  2\alpha+\frac{1-\alpha}{2}(m-2)<\rho_{\alpha}(SK_{2,\frac{m-1}{2}})
\end{aligned}
$$
for $m\geq 9$ and $\alpha \in [\frac{1}{2}, 1)$.}}
\par{\textbf{Case 3.2.} $d(w)=\frac{m-1}{2}$.}
\par{By Lemma \ref{lem3}, we have $e(N(w), V(G)\backslash N(w))\geq m-1$. Combining this with Fact 1, we have $e(N(w),V(G)\backslash N(w))= m-1$, that is, $e(V(G)\backslash N[w])=1$. It follows that $d(w_i)=2$ for $w_i \in N(w)$ by Lemmas \ref{lem3} and \ref{lem4}. Let $e(V(G)\backslash N[w])={v_1v_2}$. Then $V(G)\backslash N[w]=\{v_1,v_2\}$. Otherwise $G-w$ is disconnected, a contridiction. This implies that $G=G(a,b)$, where $1\leq a\leq b$ and $a+b=\frac{m-1}{2}$ (see Fig.\;\ref{fig1}). We assume $N_{N(w)}(v_1)=\{w_1,w_2,\ldots,w_a\}$ and $N_{N(w)}(v_2)=\{w_{a+1},w_{a+2},\ldots,w_{a+b}\}$. It is easy to see that $G(1,\frac{m-3}{2})\cong SK_{2,\frac{m-1}{2}}$. Let $X=\left(x_1, x_2, \ldots, x_n\right)^{T}$ be the $\alpha$-Perron vector of $G$ corresponding to $\rho_\alpha(G)$. Without loss of generality, we assume $x_{v_1}\leq x_{v_2}$. Then $G(1,\frac{m-3}{2})=G(a,b)-\{v_1w_i:i=1,2,\ldots,a-1\}+\{v_2w_i:i=1,2,\ldots,a-1\}$. By Lemma \ref{lem7}, we have
$$\rho_\alpha(G)=\rho_\alpha(G(a,b))\leq \rho_\alpha(G(1,\frac{m-3}{2}))=\rho_{\alpha}{(SK_{2,\frac{m-1}{2}})}$$
for $m\geq 9$ and $\alpha \in [\frac{1}{2}, 1)$, with equality if and only if $a=1$, that is, $G\cong SK_{2,\frac{m-1}{2}}$. $\qedsymbol$}
\par{Combining Facts 2 and 3, we have $\rho_{\alpha}(G) \leq \rho_{\alpha}(SK_{2,\frac{m-1}{2}})$
for $m\geq 9$ and $\alpha \in [\frac{1}{2}, 1)$, with equality if and only if $G\cong SK_{2,\frac{m-1}{2}}$. We complete the proof of (ii). $\qedsymbol$}



\begin{thebibliography}{99}
\bibitem{BZ} A. Berman, X.-D. Zhang, On the spectral radius of graphs with cut vertices, J. Combin. Theory Ser. B 83 (2001) 233--240.
\bibitem{B} B. Bollob\'{a}s, Extremal Graph Theory, Academic Press, London, New York, 1978.
\bibitem{BM} J.A. Bondy, U.S.R. Murty, Graph theory, Springer, New York, 2008.
\bibitem{BH} R.A. Brualdi, A.J. Hoffman, On the spectral radius of (0,1)-matrices, Linear Algebra Appl. 65 (1985) 133--146.
\bibitem{BS}R.A. Brualdi, E.S. Solheid, On the spectral radius of complementary acyclic matrices of zeros and ones, SIAM J. Algebra. Discrete Methods 7 (1986) 265--272.
\bibitem{CZ} M.Z. Chen, X.-D. Zhang, Some new results and problems in spectral extremal graph theory (in Chinese) J. Anhui Univ. Nat. Sci. 42 (2018) 12--25.
\bibitem{CG} X.D. Chen, L.T. Guo, On minimally 2-(edge)-connected graphs with extremal spectral radius, Discrete Math.  342 (2019) 2092--2099.
\bibitem{CR} D. Cvetkovi\'{c} and P. Rowlinson, The largest eigenvalue of a graph: A survey, Linear Multilinear Algebra 28 (1990) 3--33.
\bibitem{D} G.A. Dirac, Minimally 2-connected graphs, J. Reine Angew. Math. 228 (1976) 204--216.
\bibitem{FGL} D.D. Fan, S. Goryainov, H.Q. Lin, On the (signless Laplacian) spectral radius of minimally $k$-(edge)-connected graphs for small $k$, Discrete Appl. Math. 305 (2021) 154--163.
\bibitem{FW} Z.M. Feng, W. Wei, On the $A_\alpha$-spectral radius of graphs with given size and diameter, Linear Algebra Appl. 650 (2022) 132--149.
\bibitem{GS} J.M. Guo, J.Y. Shao, On the spectral radius of trees with fixed diameter, Linear Algebra Appl. 413 (2006) 131--147.
\bibitem{GZ2} S.-G. Guo, R. Zhang, Sharp upper bounds on the $Q$-index of (minimally) 2-connected graphs with given size, Discrete Appl. Math. 320 (2022) 408--415.
\bibitem{GZ1} S.-G. Guo, R. Zhang, The sharp upper bounds on the $A_\alpha$-spectral radius of $C_4$-free graphs and Halin graphs, Graphs Combin. 38 (2022) 19.
\bibitem{LLF} Y.T. Li, W.J. Liu, L.H. Feng, A survey on spectral conditions for some extremal graph problems, Adv. Math. 51 (2) (2022) 193--258.
\bibitem{LHX}H.Q. Lin, X. Huang, J. Xue, A note on the $A_\alpha$-spectral radius of graphs, Linear Algebra Appl. 557 (2018) 430--437.
\bibitem{LLT} H.Q. Liu, M. Lu, F. Tian, On the spectral radius of unicyclic graphs with fixed diameter, Linear Algebra Appl. 420 (2007) 449--457.
\bibitem{LBW} X.X. Liu, H.J. Broersma, L.G. Wang, On a conjecture of Nikiforov involving a spectral radius condition for a graph to contain all trees, Discrete Math. 345 (2022) 113112.
\bibitem{LMH} Z.Z. Lou, G. Min, Q.X. Huang, On the spectral radius of minimally 2-(edge)-connected graphs with given size, https://arxiv.org/abs/2206.07872.
\bibitem{N2011} V. Nikiforov, Some new results in extremal graph theory, in: Surveys in Combinatories 2011, London Math. Soc. Lecture Note Ser. 392 (2011) 141--181.
\bibitem{N2017} V. Nikiforov, Merging the $A$- and $Q$-spectral theories, Appl. Anal. Discrete Math. 11 (2017) 81--107.
\bibitem{NR}V. Nikiforov, O. Rojo, On the $\alpha$-index of graphs with pendent paths, Linear Algebra Appl. 550 (2018) 87--104.
\bibitem{P} M. Plummer, On minimal blocks, Trans. Amer. Math. Soc. 134 (1968) 85--94.
\bibitem{R}P. Rowlinson, On the maximal index of graphs with a prescribed number of edges, Linear Algebra Appl. 110 (1988) 43--53.
\bibitem{Stanic} Z. Stani\'{c}, Inequalities for Graph Eigenvalues, Cambridge Unoversity Press, New York, 2015.
\bibitem{Stanley}P.R. Stanley, A bound on the spectral radius of graphs with $e$ edges, Linear Algebra Appl. 87 (1987) 267--269.
\bibitem{Ste} D. Stevanovi\'{c}, Spectral Radius of Graphs, Academic Press, New York, 2015.
\bibitem{TCC} G.X. Tian, Y.X. Chen, S.Y. Cui, The extremal $\alpha$-index of graphs with no 4-cycle and 5-cycle, Linear Algebra Appl. 619 (2021) 160--175.
\bibitem{XLL} J. Xue, H.Q. Lin, S.T. Liu, J.L. Shu, On the $A_\alpha$-spectral radius of a graph, Linear Algebra Appl. 550 (2018) 105--120.
\bibitem{ZLS} M.Q. Zhai, H.Q. Lin, J.L. Shu, Spectral extrema of graphs with fixed size: cycles and complete
bipartite graphs, Electron. J. Comb. 95 (2021) 103322.
\bibitem{ZXL}M.Q. Zhai, J. Xue, Z.Z. Lou, The signless Laplacian spectral radius of graphs with a prescribed number of edges, Linear Algebra Appl. 603 (2020) 154--165.
\bibitem{ZL} H.H. Zhang, S.C. Li, On the Laplacian spectral radius of bipartite graphs with fixed order and size, Discrete Appl. Math. 229 (2017) 139--147.
\bibitem{ZG} R. Zhang, S.-G. Guo, Maxima of the Laplacian spectral radius of (minimally) 2-connected graphs with fixed size, Linear Algebra Appl. 651 (2022) 390--406.
\end{thebibliography}
\end{document}